\documentclass{amsart}
\usepackage{amsmath, amssymb, amsthm, graphicx}

\theoremstyle{plain}
\newtheorem{prop}{Proposition}

\newtheorem{conj}[prop]{Conjecture}
\newtheorem{lemm}[prop]{Lemma}
\newtheorem{ques}[prop]{Question}

\theoremstyle{definition}
\newtheorem*{defi}{Definition}

\newtheorem*{rema}{Remark}

\numberwithin{equation}{section}

\def\Reff#1; #2; #3; #4; #5; #6; #7\par{%
\bibitem{#1} #2, {\it #3}, #4 {\bf #5} (#6) #7}

\def\Ref#1; #2; #3; #4\par{%
\bibitem{#1} #2, {\it #3}, #4}

\renewcommand{\aa}{\mathtt{a}}
\renewcommand{\AA}{\mathtt{A}}
\newcommand{\bb}{\mathtt{b}}
\newcommand{\BB}{\mathtt{B}}
\newcommand{\cc}{\mathtt{c}}
\newcommand{\CC}{\mathtt{C}}
\def\cl#1{\overline{\vrule height5pt width0pt #1}}

\let\D=\Delta
\newcommand{\del}{\mathrm{del}}
\newcommand{\DL}{D_{\scriptscriptstyle\!L}}
\newcommand{\DR}{D_{\scriptscriptstyle\!R}}
\let\e=\varepsilon
\newcommand{\eg}{{\it e.g.}}
\newcommand{\eq}{\leftrightarrow}

\newcommand{\eqred}{\leftrightarrow^{\!\scriptscriptstyle{r}}}
\newcommand{\eqstr}{\leftrightarrow^{\!\scriptscriptstyle{s}}}
\let\ge=\geqslant
\newcommand{\ie}{{\it i.e.}}
\newcommand{\inv}{^{-1}}
\newcommand{\Inv}{{}^{-1}}
\let\le=\leqslant
\newcommand{\NL}{N_{\scriptscriptstyle\!L}}
\newcommand{\NR}{N_{\scriptscriptstyle\!R}}
\let\pp=\ldots
\newcommand{\pmo}{^{\pm 1}}
\newcommand{\resp}{{\it resp.~}}

\newcommand{\rvr}{\curvearrowright}
\newcommand{\rvl}{\mathrel{\raisebox{5pt}{\rotatebox{180}{$\curvearrowleft$}}}}
\newcommand{\rvlred}{\rvl^{\!\scriptscriptstyle{r}}}
\newcommand{\rvrred}{\rvr^{\!\scriptscriptstyle{r}}}
\newcommand{\rvlstr}{\rvl^{\!\scriptscriptstyle{s}}}
\newcommand{\rvrstr}{\rvr^{\!\scriptscriptstyle{s}}}
\let\s=\sigma
\def\ss#1{\sigma_{#1}^{\phantom 1}}
\def\sss#1{\sigma_{#1}^{-1}}
\def\ssss#1{\sigma_{#1}^{\pm1}}
\let\u=\underline

\begin{document}

\title{On word reversing in braid groups}

\author{Patrick Dehornoy}
\address{Laboratoire de Math\'ematiques Nicolas
Oresme, UMR 6139 CNRS, Universit\'e de Caen, 14032
Caen, France}
\email{\tt dehornoy@math.unicaen.fr}

\author{Bert Wiest}
\address{IRMAR, UMR 6625 CNRS, Universit\'e de
Rennes 1,  Campus Beaulieu, 35042 Rennes, France}
\email{\tt bertold.wiest@math.univ-rennes1.fr}

\begin{abstract}
It has been conjectured that in a braid group, or more
generally in a Garside group, applying any
sequence of monotone equivalences and word
reversings can increase the length of a word by
at most a linear factor depending on the group
presentation only. We give a counter-example to this
conjecture, but, on the other hand, we establish
length upper bounds for the case when only right
reversing is involved. We also state a new
conjecture which would, like the above one, imply
that the space complexity of the handle reduction
algorithm is linear.
\end{abstract}



\maketitle

This paper was motivated by attempts to estimate
the complexity of the handle reduction algorithm in
braid groups \cite{Dfo}, {\it via} a detailed study of
\emph{word reversings}.

Word reversing is a general combinatorial method for
investigating monoids and groups specified by
explicit presentations \cite{Dff, Dgc, Dgp}. In good
cases, typically in the case of braid
groups~\cite{Dff} and, more generally, Garside groups
\cite{Dgk}, it provides algorithmic solutions to the
word problem, as well as an efficient way for proving
properties such as cancellativity or existence of
least common multiples in the monoid or quadratic
isoperimetric inequalities in the group.

However, many natural questions about word reversing
remain open, even in the basic case of the standard
presentation of Artin's braid group~$B_n$. There are 
two types of word reversing, namely the left and
the right one. In the case of~$B_n$ and, more
generally, in the case of Artin--Tits groups of
finite Coxeter type, Garside's theory implies
that every sequence of right reversings must
terminate, and it gives an upper bound on the
length of the final word thus obtained; however,
it says nothing about the length  of the
\emph{intermediate} words and about many related 
questions. Also, very little is known about what
happens when both the left and the right types
are used in one reversing sequence. In
particular, we raised

\begin{ques} \cite{Dgb} \label{Q:Main}
Does there exist a constant~$C_n$ such that
the length of every freely reduced braid word
obtained from a length~$\ell$ word by using left
and right reversing plus monotone
equivalence---precise definitions are given
below---is bounded above by~$C_n\ell$?
\end{ques}

A positive answer would have implied a linear
upper bound on the space complexity of the
handle reduction algorithm in braid groups, and
indeed a positive answer was carelessly  proposed
as a conjecture in~\cite{Dgr}. The aim of this 
paper is, on the one hand, to answer
Question~\ref{Q:Main}  in the negative, by proving

\begin{prop} \label{P:LRM}
Let $w$ be the $4$~strand braid word~$\sss2
\ss1 \ss3 \ss2$. Then arbitrarily long freely
reduced words can be obtained from~$w$ using left
and right reversing and monotone equivalence.
\end{prop}

In fact, the result of Proposition~\ref{P:LRM} can 
even be strengthened by requiring that all
involved  words contain no commuting pattern
like~$\ss1 \ss3
\sss1$. 

On the other hand, we shall establish some positive 
results, namely:

\begin{prop} \label{P:Bound}
Let $w$ be an $n$~strand braid word of
length~$\ell$.\\
(i) Every word obtained from	$w$  using right
reversing has length at most~$C_n \ell$, with $C_n =
\frac12 3^n$.\\
(ii) Every positive--negative word obtained from~$w$
using right reversing and monotone equivalence has
length at most~$C'_n \ell$, with $C'_n = \frac12
n(n-1)-1$.\\
(iii) Every word obtained from~$w$ using right
reversing and monotone equivalence has length at
most~$2^{C''_n \ell}$, with $C''_n = \frac12 n(n-1)$.
\end{prop}

The upper bounds of
Proposition~\ref{P:Bound}$(i)$
and~\ref{P:Bound}$(iii)$ are certainly not optimal,
but they seem to be the first ones in this direction.
As for Proposition~\ref{P:Bound}$(ii)$, we notice
that $C'_n$ has to grow at least linearly with~$n$,
as right reversing the word
$(\ss1 \ss3
\pp
\ss{2\ell-1})\inv (\ss2 \ss4 \pp \ss{2\ell})$ leads
to a positive--negative word of length~$O(\ell^2)$.

At the end of the paper we shall propose an 
alternative conjecture which does appear to be true, 
and which would still imply a linear bound on
the space complexity of the handle reduction
algorithm.

Before giving the technical definitions, we explain
in some more detail the connection of our results
with the handle reduction algorithm~\cite{Dfo} and
$\s$-definite forms of braids. It is known that
every braid word is equivalent {\it modulo} the
braid relations to a $\ss1$-definite word, \ie~a 
word in which at least one of the letters~$\ss1,
\sss1$ does not occur. This fact is one of the two 
key points in the construction of a canonical 
ordering on braids \cite{Dgr}. Handle reduction is 
a combinatorial method that solves the isotopy problem 
of braids and produces $\ss1$-definite forms. Although
extremely efficient in practice, the method remains
partly mysterious and its exact complexity is
unknown: the only upper bound proved so far is
exponential, very far from statistical evidence.

Even more frustrating is the lack of control on
the length of the words appearing in the process:
the only proved result is an exponential upper
bound, while all experiments indicate that their length
is bounded by $C_n \ell$, where $\ell$ is the length
of the input braid word, and $C_n$ in a constant
which appears to be growing linearly with the number
of strands $n$---for four strands, the choice
$C_4=2$ seems sufficient, and as the example of
the words
$$
\ss1 \sigma_2^{-2} \sigma_3^2 \sigma_4^{-2} 
\ldots
\sigma_{n-3}^{2\epsilon}
\sigma_{n-2}^{-2\epsilon}
\sigma_{n-1}^{2\epsilon} \cdot
\sigma_{n-2}^{2\epsilon}
\sigma_{n-3}^{-2\epsilon} \sigma_{n-4}^{2\epsilon} 
\ldots
\sigma_3^{-2} \sigma_2^2 \sigma_1^{-1}$$
(with $\epsilon = \pm 1$ according to the
parity of~$n$) demonstrates,
$C_n$ needs to grow
\emph{at least}  linearly with $n$.
Now, since handle reduction is a compound
of reversing and monotone equivalence, an affirmative
answer to Question~\ref{Q:Main} would have given
the expected linear bound for the length of the
words appearing in handle reduction. As a
corollary, it would have shown that, for fixed $n$,
every braid word of length~$\ell$ is equivalent to a
$\ss1$-definite word of
length~$O(\ell)$. Let us mention that
the latter statement has been proved recently
in~\cite{DyW} using a deep result about train
tracks \cite{Ham}. The current results leave the
questions about handle reduction open. However,
handle reduction is in fact a compound of a
more restricted set of operations, namely
reversings and commutation relations, so we would be
satisfied if the length of words remained bounded
under iterated  applications of these two
operations---this is exactly the modified
conjecture stated at the end of
the paper.


\section{Word reversing}

The standard presentation of Artin's $n$~strand braid
group~$B_n$ is
$$\langle\ss1, \pp, \ss{n-1} \,;\,
\ss i \ss j = \ss j \ss i \text{~for $\vert i -
j\vert \ge 2$},
\   \ss i \ss j \ss i  = \ss j \ss i \ss j
\text{~for
$\vert i - j\vert = 1$} \rangle.$$
We denote by~$B_n^+$ the monoid with the above
presentation. An $n$~strand {\it braid word} is a word
on the $2n-2$ letters $\ssss1, \pp, \ssss{n-1}$. We
say that a braid word~$w$ is {\it positive} (\resp
{\it negative}) if no letter~$\sss i$ (\resp $\ss i$)
occurs in~$w$. We say that $w$ is {\it
positive--negative} if $w$ consists of positive
letters followed by negative letters, \ie, if
$w$ can be expressed as~$uv\inv$ with $u, v$
positive.

The operations we study here are the
following transformations on braid words:

\begin{defi}
Let $w, w'$ be braid words.

$(i)$ We say that $w$ is {\it right reversible}
to~$w'$, denoted~$w \rvr w'$, if one can transform~$w$ 
to~$w'$ by (iteratively) replacing some 
subword~$\sss i\ss j $ with $\ss j \sss i$ (case 
$\vert i - j\vert \ge 2$), or with 
$\ss j \ss i \sss j \sss i$ (case 
$\vert i - j\vert = 1$), or with~$\e$ (the empty
word, case $i = j$).

$(ii)$ Symmetrically, we say that $w$ is {\it left
reversible} to~$w'$, denoted~$w \rvl w'$, if $w'$ 
is obtained by (iteratively) 
replacing some subword~$\ss i \sss j$ with 
$\sss j \ss i$ (case $\vert i - j\vert \ge 2$), 
with $\sss j \sss i \ss j \ss i$ 
(case $\vert i - j\vert = 1$), or with~$\e$ 
(case $i = j$).

$(iii)$ We say that $w$ and $w'$ are {\it
monotonously equivalent}, denoted~$w \eq w'$, if
$w'$ is obtained from~$w$ by (iteratively)
replacing some subword~$(\ss i \ss j)\pmo$ with
$(\ss j \ss i)\pmo$ (case $\vert i - j\vert \ge
2$), or some subword~$(\ss i \ss j \ss i)\pmo$
with $(\ss j\ss i \ss j)\pmo$ (case $\vert i -
j\vert = 1$).
\end{defi}

It is clear that reversing and monotone equivalence
transforms a braid word into an equivalent word,
\ie, one that represents the same element of the
braid group. Observe that the above transformations
never introduce trivial pairs of the form $\ss i \sss
i$ or $\sss i \ss i$. So, typically, for a braid
word~$w$ to be reversible to the empty word~$\e$
is {\it a priori} a stronger condition than just
being equivalent to~$\e$, as one is allowed to
introduce no $\ss i \sss i$ or $\sss i \ss i$ in order
to transform~$w$ into~$\e$.

Clearly, the words that cannot be transformed using
right reversing are the positive--negative words. The
key result about braid word reversing is as follows:

\begin{prop} \label{P:Basic} \cite{Dff}
Let $w$ be an $n$~strand braid word of
length~$\ell$. Then there exists a unique
positive--negative word~$w'$ such that $w$ is right
reversible to~$w'$. Moreover, the length of~$w'$ is 
at most~$C_n \ell$, with $C_n = \frac12 n(n-1)-1$.
\end{prop}

Proposition~\ref{P:Basic} is a consequence of
Garside's result that common right multiples
exist in braid monoids~\cite{Gar} and of
general properties of word reversing~\cite{Dgp}
guaranteeing that, for all positive words~$u,
v$, the existence of positive words~$u_1, v_1$
satisfying $u\inv v \rvr v_1 u_1\inv$ is
equivalent to the existence of a common right
multiple for the elements represented
by~$u$ and~$v$. In the current paper, we shall
only use the following result:

\begin{lemm} \label{L:Lcm} \cite{Dff}
Assume that $u, u'$ are equivalent positive
braid words and, similarly, that $v, v'$ are
equivalent positive braid words. Let $u_1, v_1,
u'_1, v'_1$ be the positive words satisfying 
$u\inv v \rvr v_1 u_1\inv$ and $u'\Inv v' \rvr
v'_1 u'_1\Inv$. Then $u_1$ and~$u'_1$ are
equivalent, and so are~$v_1$ and~$v'_1$. 
\end{lemm}

\begin{rema}
The previous results imply that
right reversing solves the word problem of the
braid monoid and of the braid group, in one and
two passes respectively. Indeed Lemma~\ref{L:Lcm}
implies that two positive braid words $u, v$
represent the same element of the braid monoid if
and only if $u\inv v$ is right reversible to the
empty word, and that an arbitrary braid word~$w$
represents~$1$ in the braid group if and only if
it is right reversible to some positive--negative
word~$v u\inv$ such that $u\inv v$ is right
reversible to the empty word. The last step is
equivalent to $v u\inv$ being left reversible to
the empty word. So a braid word~$w$
represents~$1$ if and only if the empty word can
be obtained from~$w$ using left and right
reversing.
\end{rema}


\section{Counterexamples}

Proposition~\ref{P:Basic} says nothing about
the words one obtains using both left and right
reversing. The trivial example
\begin{equation} \label{E:TrivialCE}
\sss1 \ss2 \rvr \ss2 \ss1 \sss2 \sss1
\rvl \ss2 \sss2 \sss1 \ss2 \ss1 \sss1
\end{equation}
shows that, starting from~$\sss1 \ss2$, we can
produce words of arbitrary length using left and
right reversing, since the initial word is a proper
factor of the final word. Hence, whenever
both left and reversing are involved, restricting to
freely reduced words, \ie, containing no
pattern~$\ss i \sss i$ or~$\sss i \ss i$, is a
minimal requirement if one is to expect bounded
length.

\begin{defi}
We define {\it reduced right reversing},
denoted~$\rvrred$, to be the variant of right
reversing in which a free reduction is performed
after each reversing step. Reduced left reverving
and monotone equivalence are defined similarly.
\end{defi}

Question~\ref{Q:Main} asks in particular whether the
words obtained from a given word using reduced
reversing and monotone equivalence have a bounded
length. We now establish Proposition~\ref{P:LRM},
which provides a negative answer. To improve
readability, we adopt a convention of~\cite{Eps},
using $\aa,
\bb,
\pp$ for~$\ss1,
\ss2, \pp$, and $\AA, \BB, \pp$ for
$\sss1, \sss2, \pp$. For instance,
\eqref{E:TrivialCE} becomes
$\AA \bb \rvr \bb \aa \BB \AA
\rvl \bb \BB \AA \bb \aa \AA$.

\begin{figure} [htb]
\begin{picture}(126,75)(0, 0)
\put(0,0){\includegraphics{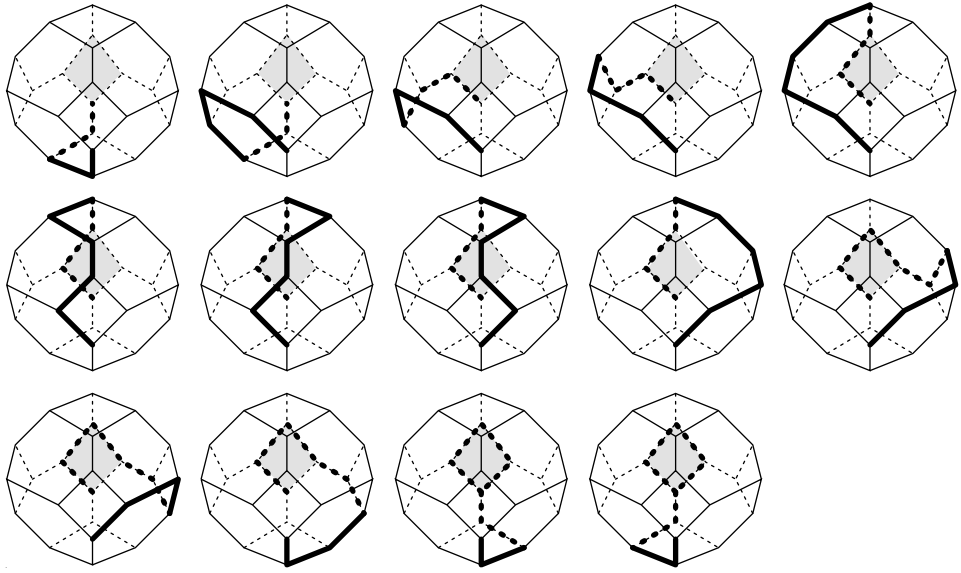}}
\end{picture}
\caption{\smaller Generating arbitrarily long
words from $\sss2 \s_1 \s_3 \s_2$ using reversing
and monotone equivalence; all words are traced on
the fragment of the Cayley graph corresponding to
the divisors of~$\D_4$, \ie, on the
$4$-permutohedron, which, topologically, is a
sphere; the initial path is pushed around the sphere
so as to make a loop around the grey facet on the
rear; each other facet is crossed once.}
\label{F:Spheres}
\end{figure}

\begin{proof} [Proof of Proposition~\ref{P:LRM}]
(Figure~\ref{F:Spheres}) We find (the underlined
subwords are those we transform):
\begin{align*}
\u{\BB \aa} \cc \bb
& \rvrred \aa \bb \AA \u{\BB \cc \bb}
  \rvrred \aa \bb \u{\AA \cc} \bb \CC
  \rvrred \aa \bb \cc \u{\AA \bb} \CC
  \rvrred \aa \u{\bb \cc \bb} \aa \BB \AA \CC \\
&  \eqred \u{\aa \cc} \bb \cc \aa \BB \AA \CC  
  \eqred \cc \aa \bb \u{\cc \aa} \BB \AA \CC 
  \eqred \cc \u{\aa \bb \aa} \cc \BB \AA \CC 
 \eqred \cc \bb \aa \u{\bb \cc \BB} \AA \CC \\
&  \rvlred \cc \bb \u{\aa \CC} \bb \cc \AA \CC
  \rvlred \u{\cc \bb \CC} \aa \bb \cc \AA \CC
  \rvlred \BB \cc \u{\bb \aa \bb} \cc \AA \CC
  \eqred \BB \u{\cc \aa} \bb \aa \cc \AA \CC
  \eqred \BB \aa \cc \bb \aa \cc \AA \CC,
\end{align*}
and, inductively, $\BB \aa \cc \bb$ transforms
into $\BB \aa \cc \bb (\aa \cc \AA \CC)^k$ for
each~$k$ as the the words above never finish with
the letter~$\AA$.
\end{proof}

Note that in the previous counter-example not
only the final words, but even all intermediate words
are freely reduced. Now we see that these words
still involve the commuting pattern
$\aa\cc\AA\CC$, \ie, $\ss1\ss3\sss1\sss3$. We
shall show now that even such semi-trivial patterns 
can be avoided.

\begin{defi}
We say that a braid word is {\it strongly reduced} if
it is freely reduced and, in addition, contains no
subword of the form $\s_i^e \s_j^d \s_i^{-e}$ with
$e, d = \pm 1$ and $\vert i - j\vert \ge 2$. We
define {\it strongly reduced} right reversing to be
the variant in which a full reduction is performed
after each reversing step.
\end{defi}

In the above definition, strongly reducing a word means
iteratively replacing each subword of the form
$\s_i^e \s_j^d \s_i^{-e}$ with the corresponding
letter~$\s_j^d$. This is easily seen to lead in
finitely many steps to a strongly reduced word. The
latter need not be unique, but the various words so
obtained are equivalent via commutation relations.

\begin{prop} \label{P:LRM+}
Starting from~$\sss1 \sss3 \sss2 \ss1 \ss2 \ss3 \ss2
\sss3 \sss2 \sss1 \ss2 \sss3$, one can derive using
strongly reduced left and right reversing and monotone
equivalence arbitrary long (strongly reduced) words.
\end{prop}

\begin{proof}
Using~$\rvrstr$, $\rvlstr$ and $\eqstr$ for the
strongly reduced versions of~$\rvr$, $\rvl$,
and~$\eq$, we find

\hspace{10mm} $\AA \CC \BB \aa \bb \cc \bb \CC 
\u{\BB \AA \bb} \CC  (\BB \aa \bb \cc \BB \CC \BB 
\cc \bb \aa \BB \AA)^k$

\hspace{28mm} $\rvrstr \AA \CC \BB \aa \bb \cc \bb 
\u{\CC \aa} \BB \AA \CC (\BB \aa \bb \cc \BB \CC 
\BB \cc \bb \aa \BB \AA)^k$

\hspace{28mm} $\rvrstr \AA \CC \BB \aa \bb \cc \bb
\aa \CC \BB \u{\AA \CC} (\BB \aa \bb \cc \BB \CC 
\BB \cc \bb \aa \BB \AA)^k $

\hspace{28mm} $\eqstr \AA \CC \BB \aa \bb \cc \bb
\aa \u{\CC \BB \CC} \AA (\BB \aa \bb \cc \BB \CC 
\BB \cc \bb \aa \BB \AA)^k$

\hspace{28mm} $\eqstr \AA \CC \BB \aa \u{\bb \cc \bb} 
\aa \BB \CC \BB \AA (\BB \aa \bb \cc \BB \CC \BB 
\cc \bb \aa \BB \AA)^k$

\hspace{28mm} $\eqstr \AA \CC \BB \aa \cc \bb \cc 
\u{\aa \BB} \CC \BB \AA (\BB \aa \bb \cc \BB \CC \BB 
\cc \bb \aa \BB \AA)^k$

\hspace{28mm} $\rvlstr \AA \CC \BB \aa \cc \bb \cc \BB
\AA \bb \u{\aa \CC} \BB \AA (\BB \aa \bb \cc \BB 
\CC \BB \cc \bb \aa \BB \AA)^k$

\hspace{28mm} $\rvlstr \AA \CC \BB \aa \cc \bb \cc \BB
\AA \u{\bb \CC} \aa \BB \AA (\BB \aa \bb \cc \BB \CC 
\BB \cc \bb \aa \BB \AA)^k$

\hspace{28mm} $\rvlstr \AA \CC \BB \u{\aa \cc} \bb
\cc 
\BB \AA \CC \BB \cc \bb \aa \BB \AA (\BB \aa \bb \cc 
\BB \CC \BB \cc \bb \aa \BB \AA)^k$

\hspace{28mm} $\eqstr \AA \u{\CC \BB} \cc \aa \bb \cc 
\BB \AA \CC \BB \cc \bb \aa \BB \AA (\BB \aa \bb \cc 
\BB \CC \BB \cc \bb \aa \BB \AA)^k$

\hspace{28mm} $\eqstr \AA \bb \CC \u{\BB \aa \bb} \cc 
\BB \AA \CC \BB \cc \bb \aa \BB \AA (\BB \aa \bb \cc 
\BB \CC \BB \cc \bb \aa \BB \AA)^k$

\hspace{28mm} $\rvrstr \AA \bb \CC \aa \bb \u{\AA \cc} 
\BB \AA \CC \BB \cc \bb \aa \BB \AA (\BB \aa \bb \cc 
\BB \CC \BB \cc \bb \aa \BB \AA)^k$

\hspace{28mm} $\rvrstr \AA \bb \CC \aa \bb \cc \u{\AA 
\BB \AA} \CC \BB \cc \bb \aa \BB \AA (\BB \aa \bb \cc 
\BB \CC \BB \cc \bb \aa \BB \AA)^k$

\hspace{28mm} $\eqstr \AA \bb \CC \aa \u{\bb \cc \BB}
\AA \BB \CC \BB \cc \bb \aa \BB \AA (\BB \aa \bb \cc
\BB \CC \BB \cc \bb \aa \BB \AA)^k$

\hspace{28mm} $\rvlstr \AA \bb \CC \aa \u{\CC \bb
\cc} 
\AA \BB \CC \BB \cc \bb \aa \BB \AA (\BB \aa \bb \cc 
\BB \CC \BB \cc \bb \aa \BB \AA)^k$

\hspace{28mm} $\rvrstr \AA \u{\bb \CC} \aa \bb \cc \BB 
\AA \BB \CC \BB \cc \bb \aa \BB \AA (\BB \aa \bb \cc 
\BB \CC \BB \cc \bb \aa \BB \AA)^k$

\hspace{28mm} $\rvlstr \AA \CC \BB \cc \u{\bb \aa \bb}
\cc
  \BB \AA \BB \CC \BB \cc \bb \aa \BB \AA (\BB \aa \bb
\cc \BB \CC \BB \cc \bb \aa \BB \AA)^k$

\hspace{28mm} $\eqstr \u{\AA \CC} \BB \cc \aa \bb \aa
\cc
\BB \AA \BB \CC \BB \cc \bb \aa \BB \AA (\BB \aa \bb
\cc \BB \CC \BB \cc \bb \aa \BB \AA)^k$

\hspace{28mm} $\eqstr \CC \AA \BB \cc \aa \bb \aa
\u{\cc
\BB} \AA \BB \CC \BB \cc \bb \aa \BB \AA (\BB \aa \bb
\cc \BB \CC \BB \cc \bb \aa \BB \AA)^k$

\hspace{28mm} $\rvlstr \CC \AA \BB \u{\cc \aa} \bb
\aa 
\BB \CC \bb \cc \AA \BB \CC \BB \cc \bb \aa \BB \AA 
(\BB \aa \bb \cc \BB \CC \BB \cc \bb \aa \BB \AA)^k$

\hspace{28mm} $\eqstr \CC \u{\AA \BB \aa} \cc \bb \aa 
\BB \CC \bb \cc \AA \BB \CC \BB \cc \bb \aa \BB \AA 
(\BB \aa \bb \cc \BB \CC \BB \cc \bb \aa \BB \AA)^k$

\hspace{28mm} $\rvrstr \CC \bb \AA \u{\BB \cc \bb} \aa 
\BB \CC \bb \cc \AA \BB \CC \BB \cc \bb \aa \BB \AA 
(\BB \aa \bb \cc \BB \CC \BB \cc \bb \aa \BB \AA)^k$

\hspace{28mm} $\rvrstr \CC \bb \AA \cc \bb \CC \u{\aa 
\BB} \CC \bb \cc \AA \BB \CC \BB \cc \bb \aa \BB \AA 
(\BB \aa \bb \cc \BB \CC \BB \cc \bb \aa \BB \AA)^k$

\hspace{28mm} $\rvlstr \CC \u{\bb \AA} \cc \bb \CC
\BB 
\AA \bb \aa \CC \bb \cc \AA \BB \CC \BB \cc \bb \aa 
\BB \AA (\BB \aa \bb \cc \BB \CC \BB \cc \bb \aa \BB
\AA)^k$

\hspace{28mm} $\rvlstr \u{\CC \AA} \BB \aa \bb \cc
\bb 
\CC \BB \AA \bb \aa \CC \bb \cc \AA \BB \CC \BB \cc \bb
\aa \BB \AA (\BB \aa \bb \cc \BB \CC \BB \cc \bb \aa
\BB \AA)^k$

\hspace{28mm} $\eqstr \AA \CC \BB \aa \bb \cc \bb \CC
\BB \AA \bb \u{\aa \CC} \bb \cc \AA \BB \CC \BB \cc 
\bb \aa \BB \AA (\BB \aa \bb \cc \BB \CC \BB \cc \bb 
\aa \BB \AA)^k$

\hspace{28mm} $\rvlstr \AA \CC \BB \aa \bb \cc \bb \CC
\BB \AA \bb \CC \aa \bb \u{\cc \AA} \BB \CC \BB \cc 
\bb \aa \BB \AA (\BB \aa \bb \cc \BB \CC \BB \cc \bb 
\aa \BB \AA)^k$

\hspace{28mm} $\rvlstr \AA \CC \BB \aa \bb \cc \bb \CC
\BB \AA \bb \CC \u{\aa \bb \AA} \cc \BB \CC \BB \cc 
\bb \aa \BB \AA (\BB \aa \bb \cc \BB \CC \BB \cc \bb 
\aa \BB \AA)^k$

\hspace{28mm} $\rvlstr \AA \CC \BB \aa \bb \cc \bb \CC
  \BB \AA \bb \CC \BB \aa \bb \cc \BB \CC \BB \cc \bb
\aa \BB \AA (\BB \aa \bb \cc \BB \CC \BB \cc \bb \aa
\BB \AA)^k$.\\
The latter word is $\AA \CC \BB \aa \bb \cc \bb \CC
\BB \AA \bb \CC (\BB \aa \bb \cc \BB \CC \BB \cc \bb
\aa \BB \AA)^{k+1}$.
\end{proof}


\section{Length upper bounds for right reversing}

Now we turn to positive results, and establish some
upper bounds for the length of the words that can be
constructed using reversing and monotone equivalence.
In this section, we consider the case of right
reversing alone. Proposition~\ref{P:Basic}
provides an upper bound on the length of the final,
\ie, positive--negative, word that can be
obtained from a word~$w$, but it gives no bound for
the intermediate words. This is what
Proposition~\ref{P:Bound}$(i)$ does.

In order to prove the result, we need some auxiliary
notions. First, as usual, we associate with each
$n$~strand braid word~$w$ the braid diagram
obtained by concatenating the elementary diagrams
for the successive letters of~$w$, and
the diagram for~$\ss i$ is
\begin{center}
\begin{picture}(55,9)(0, 0)
\put(9.5,0){\includegraphics{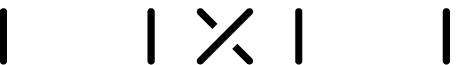}}
\put(9,7){$1$}
\put(29,7){$i$}
\put(33,7){$i\!+\!1$}
\put(54,7){$n$}
\put(0,3){$\ss i$:}
\put(16,3){$\pp$}
\put(46,3){$\pp$}
\end{picture}
\end{center}
An $n$~strand braid diagram can be seen as
the projection on $y = 0$ of a 3D-figure consisting
of $n$~non-intersecting curves.

\begin{figure} [htb]
\begin{picture}(67,33)(0, 0)
\put(0,0){\includegraphics{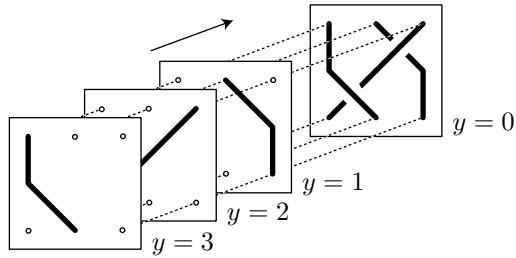}}
\put(19,0){$y= 3$}
\put(29,4){$y= 2$}
\put(39,8){$y= 1$}
\put(59,16){$y= 0$}
\end{picture}
\caption{The braid word~$\s_2 \s_1\inv$ is
layered: the strands of the associated diagram live
in parallel vertical planes}
\label{F:Simple}
\end{figure}

\begin{defi} (Figure~\ref{F:Simple})
A braid word~$w$ is said to be {\it layered} if the
associated diagram can be realized as the projection
of a 3D-figure in which each strand lives in some
vertical plane.
\end{defi}

It is well-known that, if $u$ is an $n$~strand
positive word, then $u$ is layered if and only if $u$
is {\it simple}, \ie, it represents a divisor of
Garside's fundamental braid~$\D_n$ in the
monoid~$B_n^+$.

\begin{lemm} \label{L:SimpleClosed}
(i) If $u$ and $v$ are layered positive words, then
$u\inv v$ is layered.

(ii) If $w$ is a layered word, then every word
obtained from~$w$ using reversing or monotone
equivalence is still layered.
\end{lemm}

\begin{figure} [htb]
\begin{picture}(53,24)(0, 0)
\put(0,0){\includegraphics{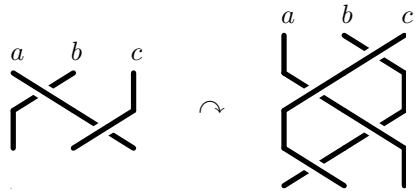}}
\put(0,17){$a$}
\put(8,17){$b$}
\put(16,17){$c$}
\put(36,22){$a$}
\put(44,22){$b$}
\put(52,22){$c$}
\put(25,10){$\rvr$}
\end{picture}
\caption{\smaller Layered words are closed under
right reversing: if the pattern $\sss1 \s_2$ occurs
in a layered word, then, necessarily, the strand~$c$
lies in the front plane, while $b$ lies in the back
plane; then the pattern $\s_2 \s_1 \sss2 \sss1$ can
be realised in the same planes without level
obstruction. The case of $\sss 2 \s_1$ is similar.}
\label{F:Closure}
\end{figure}

\begin{proof}
$(i)$ If $u$ is a positive layered word, then the
diagram of~$u$ can be realized so that the $i$-th
strand, \ie, the strand that starts at position~$i$,
lives in the plane~$y = n-i$. Thus $u\inv$ can be
realized so that the strand finishing at
position~$i$ lives in~$y = i$, and $v$ can be
realized so that the strand starting at
position~$i$ lives in the same plane. Hence the two
diagrams can be concatenated without contradicting
layeredness.

For~$(ii)$, it suffices to check that each
elementary transformation introduces no obstruction
to the hypothesis that the strands live in a
vertical plane. The case of commutation relations
is trivial. The case of right reversing is illustrated
in Figure~\ref{F:Closure}; the cases of left
reversing and monotone equivalence are similar.
\end{proof}

In the braid diagram associated with a layered
word~$w$, there is a well-defined rear strand, \ie,
the strand that lives in the plane $y = c$ with
minimal~$c$: to avoid ambiguity, we choose the
leftmost strand in that plane if there are several
ones---this makes sense as the strands living in a
given plane may not intersect.

\begin{defi}
For~$w$ a layered braid word, we denote by $\del(w)$
the braid word that encodes the diagram obtained
from the diagram of~$w$ by deleting the rear
strand.
\end{defi}

\begin{lemm} \label{L:Del}
If $w$ is a layered braid word, then $w \rvr
w'$ implies $\del(w) \rvr \del(w')$.
\end{lemm}

\begin{proof}
Once again, it is sufficient to consider the
possible elementary transformations. Up to a
translation of indices, the only non-trivial cases
are $\sss1 \ss2 \rvr \ss2 \ss1 \sss2 \sss1$
and $\sss2 \ss1 \rvr \ss1 \ss2 \sss1 \sss2$, which
both reduce to $\ss1 \rvr \ss1$ when the rear
strand is deleted.
\end{proof}

Observe that the assumption that the removed strand
is the rear (or the front) one is necessary: if we
remove the middle strand in $\sss1 \ss2 \rvr \ss2
\ss1 \sss2 \sss1$, we obtain $\e$ on the left, and
$\ss1 \sss1$ on the right. However, $\e \rvr \ss1
\sss1$ fails.

\begin{defi}
A braid word is said to be {\it elementary} if it is
a subword of a word obtained by right reversing from
a word of the form~$u\inv v$, with $u, v$ layered
positive (or, equivalently, simple) words.
\end{defi}

By Lemma~\ref{L:SimpleClosed}, every elementary word
is layered, but the converse is not true: $\ss1 \sss1$
is layered, but not elementary. Indeed, when $u\inv
v$ is right reversible to~$v_1 u_1\inv$, then the
braids represented by~$u_1$ and~$v_1$ have no common
right divisor in the braid monoid.

\begin{lemm} \label{LengthElem}
The length of an $n$~strand elementary braid word is
at most $\frac 1 2 3^n$.
\end{lemm}

\begin{figure} [htb]
\begin{picture}(75,21)(0, 0)
\put(0,0){\includegraphics{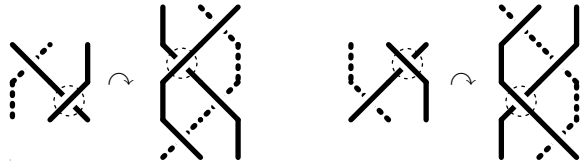}}
\put(13,10){$\rvr$}
\put(58.5,10){$\rvr$}
\end{picture}
\caption{\smaller When the rear
strand is involved in a right reversing step, it
necessarily goes through a crossing of the front
strands, from left to right.}
\label{F:Crossing}
\end{figure}

\begin{proof}
Assume that $w$ is right reversible to~$w'$ in one
step. Then the crossings in the diagram encoded
by~$w'$ are not exactly the same as the crossings in
the diagram encoded by~$w$, but we can define a
notion of inheritance: for instance, in
Figure~\ref{F:Closure}, we say that the crossing of
the strands~$a$ and~$b$ in the right figure is the
heir of the crossing of these strands in the left
figure. Then it is easy to check that, in each
case, the crossings in~$w'$ are the heirs of the
crossings of~$w$, except that two new crossings may
appear (\eg, crossings of~$b$ and~$c$ in
Figure~\ref{F:Closure}), or two crossings may vanish
(when a free reduction is performed).

Let us consider a right reversing sequence $w_0$, \pp,
$w_r$, \ie, we assume that $w_k$ is right
reversible to~$w_{k+1}$ in one step for each~$k$. 
We define the {\it total number of crossings}~$C$ in
this sequence as follows: each crossing in~$w_0$
contributes~$1$ to~$C$, and so does every new crossing
that appears in some~$w_k$, even if it subsequently
vanishes; on the other hand, the contribution to~$C$
of a crossing that is the heir of a previously
existing crossing is~$0$. So $C$ is the sum of the
number of crossings of all the terms in the sequence,
up to inheritance.

We claim that $E_n = \frac{1}{2} 3^n - n - \frac{1}{2}$
yields an upper bound for the total number of
crossings in a right reversing sequence starting with
a (layered) word of the form~$u\inv v$ with $u, v$
simple $n$~strand braid words. Then, in particular,
$E_n$ is an upper bound for the length of each braid
word occurring in such a reversing sequence, and,
therefore, for the length of every $n$~strand
elementary braid word.

For~$n = 2$, the only sequence to consider is~$(\sss1
\ss1, \e)$, so $E_2 = 2$ is indeed a valid upper bound.

Assume $n \ge 3$, and let $w$ be an $n$~strand
elementary word. By hypothesis there is a finite
sequence of words $w_0 = u\inv v$, $w_1$, \pp, $w_r =
w$ such that $u, v$ are positive layered words and
each word~$w_i$ is right reversible to~$w_{i+1}$ in
one step. By Lemma~\ref{L:SimpleClosed}, all the
words~$w_k$ are layered. By Lemma~\ref{L:Del}, the
words $\del(w_0)$,
\pp, $\del(w_r)$ also form a right reversing sequence.
Moreover, $\del(w_0)$, \ie, $\del(u\inv v)$, is a
word of the form $u_*\inv v_*$, where $u_*, v_*$ are
layered positive $n-1$~strand words. So each
word~$\del(w_k)$ is elementary, and, by induction
hypothesis, the total number of crossings (up to
inheritance) in the sequence $\del(w_0)$, \pp,
$\del(w_k)$ is bounded above by~$E_{n-1}$.

Now let us reintroduce the rear strand in the
initial word and count how many crossings it can
create in (the diagram associated with)~$w_i$.
First, in~$u$ and~$v$, which are positive, the rear
strand may cross each other strand at most once, so
it creates at most $2(n-1)$ crossings. 
Then, the reversing steps may create new
crossings between the front stands and the rear
strand. However, we claim that at most $2E_{n-1}$
such crossings can be created during the sequence of
right reversings. Indeed, during each such reversing,
the rear strand moves behind one crossing of the
remaining strands, from left to right, and in the
process it creates two new crossings
(Figure~\ref{F:Crossing}). As the total number of
crossings not involving the rear strand is at most
$E_{n-1}$, this puts the desired bound on the number
of new reversings. 
In summary, we obtain $E_n \le E_{n-1} + 2(n-1) + 2
E_{n-1}$. Now we calculate
$$3E_{n-1} + 2(n-1) = 
 \frac{1}{2} 3^n - 3n - \frac{3}{2} + 2n - 2 
 < \frac{1}{2} 3^n - n - \frac{1}{2}.
$$
This completes the proof of the lemma.
\end{proof}

We remark that the proof of lemma \ref{LengthElem} 
would go through even if we allowed not only right
reversings, but also commutation relations like
$\sigma_1 \sigma_3 \mapsto \sigma_3 \sigma_1$ or even
$\sigma_1 \sigma_3^{-1} \mapsto 
\sigma_3^{-1} \sigma_1$. 

Now it remains to decompose arbitrary words into
products of elementary words. We recall that a layered
positive word is the same as a simple word, in the
sense of Garside, namely a positive word representing
a divisor of $\Delta_n$. From now on, we shall be
dealing with positive words only, and therefore use
the word ``simple'' rather than ``layered''.

\begin{lemm} \label{L:Decomposition}
Assume $w = w_1^{\epsilon_1} \pp w_\ell^{\epsilon_\ell}$, where
$w_1, \pp, w_\ell$ are simple positive braid words,
and $\epsilon_1, \pp, \epsilon_\ell = \pm1$. Then every word~$w'$
obtained from~$w$ using right reversing can be
written as the product of at most~$\ell$ elementary
words.
\end{lemm}

\begin{proof}
First we associate with every right reversing sequence
$\widetilde{w}_0$, $\widetilde{w}_1$, \pp a planar
oriented graph whose edges are labeled by~$\ss i$'s.
This graph, which will be called a {\it reversing
diagram}, is analogous to a van Kampen diagram, and
it is constructed inductively as follows
(Figure~\ref{F:ReversingDiagram}). First we
associate with~$\widetilde{w}_0$ a path shaped like an
ascending staircase by reading~$\widetilde{w}_0$ from
left to right and iteratively appending a horizontal
right-oriented edge labeled~$\ss i$ for each
letter~$\ss i$, and a vertical down-oriented edge
labeled~$\ss i$ for each letter~$\sss i$. Assume that
the fragment corresponding to~$\widetilde{w}_0,
\pp, \widetilde{w}_{k-1}$ has been constructed and its right
side is a path labeled~$\widetilde{w}_{k-1}$. By definition,
the word~$\widetilde{w}_k$ is obtained from~$\widetilde{w}_{k-1}$ by
replacing some subword~$\sss i \ss j$ with the
unique word~$u v\inv$ such that $\ss i v = \sss j u$
is a relation of the considered presentation. The
involved subword~$\sss i \ss j$ corresponds to some
top-left oriented corner in the diagram, and we
complete the diagram and transform this corner into a
square by adding horizontal edges labelled~$u$ and
vertical edges labelled~$v$, following the scheme:
\begin{center}
\begin{picture}(56,13)(0,0)
\put(4,2.5){\includegraphics{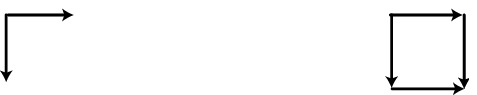}}
\put(0,6.5){$\ss i$}
\put(7,12.5){$\ss j$}
\put(40,6.5){$\ss i$}
\put(47,12.5){$\ss j$}
\put(47,1){$v$}
\put(52,6.5){$u$}
\put(14,6.5){completed into}
\end{picture}.
\end{center}

\begin{figure} [htb]
\begin{picture}(108,38)(0, 0)
\put(4,2.5){\includegraphics{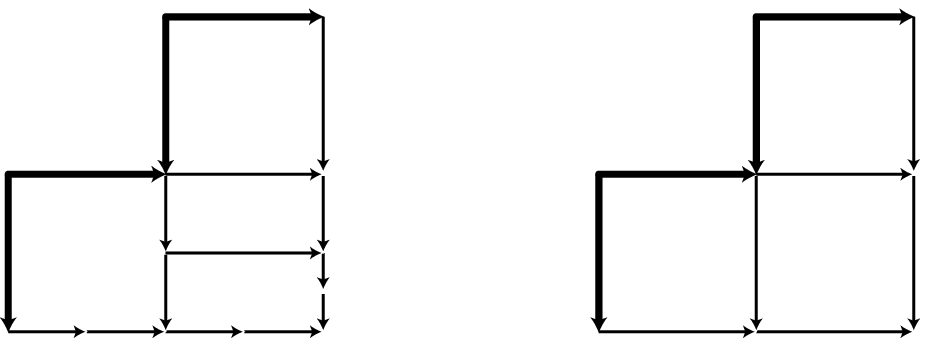}}
\put(7,0){$\ss2$}
\put(15,0){$\ss1$}
\put(23,0){$\ss3$}
\put(31,0){$\ss2$}
\put(70,0){$\ss2\ss1$}
\put(86,0){$\ss3\ss2$}
\put(28,13){$\ss3$}
\put(11,21.5){$\ss2$}
\put(28,21.5){$\ss3$}
\put(71,21.5){$\ss2$}
\put(88,21.5){$\ss3$}
\put(28,37.5){$\ss3$}
\put(88,37.5){$\ss3$}
\put(0,11){$\ss1$}
\put(22,7){$\ss2$}
\put(22,15){$\ss1$}
\put(22,27){$\ss1$}
\put(38,9){$\ss2$}
\put(38,5){$\ss3$}
\put(38,15){$\ss1$}
\put(38,27){$\ss1$}
\put(60,11){$\ss1$}
\put(82,11){$\ss1\ss2$}
\put(82,27){$\ss1$}
\put(98,11){$\ss1\ss2\ss3$}
\put(98,27){$\ss1$}
\end{picture}
\caption{\smaller The reversing diagram (left) and the
reversing grid (right) associated
with the sequence $\sss1 \s_2 \sss1 \s_3
\rvr \s_2 \s_1 \sss2 \sss1 \sss1 \s_3
\rvr \s_2 \s_1 \sss2 \sss1 \s_3 \sss1
\rvr \s_2 \s_1 \sss2 \s_3 \sss1 \sss1
\rvr \s_2 \s_1 \s_3 \s_2 \sss3 \sss2 \sss1 \sss1$:
one draws a zigzag path labelled by the initial
word, and, then, one iteratively fills the open
top-left corners using the braid relations.}
\label{F:ReversingDiagram}
\end{figure}

The next step is to observe that each right
reversing graph starting with~$w$, in particular
the maximal one, \ie, the one that finishes with
a positive--negative word, admits a rectangular 
spine, which will be called the {\it right 
reversing grid} of~$w$. Assume that $w$
contains $q$~positive letters and
$p$~negative ones.  We first assume in addition
that $w$ is a negative--positive word. We define two
sequences of simple words~$u_{i, j}, v_{i, j}$ for $0
\le i \le p$ and $0 \le j \le q$ by setting $u_{0,i}$
to be the $i$th letter in the inverse of the negative
part of~$w$ for $i < p$, and to be the empty
word~$\e$ for $i = p$, and by defining $v_{j, 0}$
to be the $j$th letter in the positive part
of~$w$ for $j < q$, and to be~$\e$ for~$j = q$. Then
we inductively define~$u_{i, j+1}, v_{i+1, j}$ by
$u_{i, j}\inv v_{i, j} \rvr v_{i+1, j} u_{i,
j+1}\inv$. We notice that these words are indeed
simple. In this way, we obtain a grid, which is a
fragment of the complete reversing diagram associated
with~$w$ (Figure~\ref{F:Grid}). If $w$ is not
negative--positive, then the construction is similar,
except that the word~$w$ need not correspond to a
top--left corner, and the top--left corner of the
rectangular grid may be missing (as in 
Figure~\ref{F:ReversingDiagram}).

\begin{figure} [htb]
\begin{picture}(66,44)(0, 0)
\put(6,2.2){\includegraphics{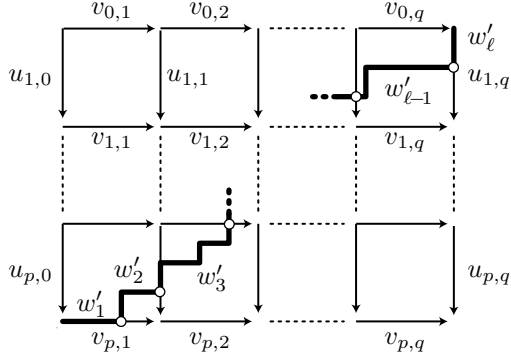}}
\put(0,9){$u_{p,0}$}
\put(0,35){$u_{1,0}$}
\put(21,35){$u_{1,1}$}
\put(61,9){$u_{p,q}$}
\put(61,35){$u_{1,q}$}
\put(11,0){$v_{p,1}$}
\put(24,0){$v_{p,2}$}
\put(50,0){$v_{p,q}$}
\put(11,44){$v_{0,1}$}
\put(24,44){$v_{0,2}$}
\put(50,44){$v_{0,q}$}
\put(11,26.7){$v_{1,1}$}
\put(24,26.7){$v_{1,2}$}
\put(50,26.7){$v_{1,q}$}
\put(9,4.5){$w'_1$}
\put(14,9){$w'_2$}
\put(25,8){$w'_3$}
\put(50,33){$w'_{\ell\!-\!1}$}
\put(61,40){$w'_{\ell}$}
\end{picture}
\caption{\smaller The reversing grid and the
decomposition of a word~$w'$ into a product of
elementary words (here delimited by white dots); in
the case of the word
$\sss1 \s_2 \sss1 \s_3$ of
Figure~\ref{F:ReversingDiagram}, we have
$u_{2, 0} = u_{1, 1}= \s_1$, $v_{1, 1} = \s_2$,
$v_{0,2} = \s_3$, and, for instance, $u_{2,2} =
\s_1 \s_2 \s_3$ and $v_{2,2} = \s_3 \s_2$.}
\label{F:Grid}
\end{figure}

The point now is that every word~$w'$ obtained
from~$w$ using right reversing labels a path
from the bottom-left corner to the top-right corner
in the right reversing diagram of~$w$. As the
reversing grid partitions this diagram into
squares, we can attribute a square of the grid to
each letter in~$w'$; we take the convention that,
when a letter~$\ss i$ corresponds to the vertical
common edge between two squares of the grid, it is
attached to the rightmost square, and that,
similarly, a horizontal edge belongs, in case of
doubt, to the bottom square. As the path labelled~$w$
contains only right-oriented edges and top-oriented
edges---actually bottom-oriented edges that are
crossed in the wrong direction---the only possibility
after a letter attached to the $(i, j)$-square is a
letter attached to $(i + e, j - d)$ with $e, d =
0$~or~$1$. It follows that the path associated
with~$w'$ visits at most $\ell$~squares, and that
$w'$ can be decomposed into a product $w'_1 \pp
w'_\ell$, where $w'_k$ consists of all letters
attached to the $k$th visited square.

It remains to see that each word~$w'_k$ is
elementary. Now an easy induction shows that, if
$w_*$ is a fragment of~$w'$ lying in the 
$(i,j)$-square, then there exist positive 
words~$u_*,v_*$ such that $u_{i, j}\inv v_{i, j}$ is 
right reversible to $u_*^e w_* v_*^d$, where $e$ 
and~$d$ are~$\pm 1$ according to the sides of the 
square through which the path enters and exits the 
square. Since the words $u_{i, j}$ and $v_{i,j}$ are
simple, the words~$u_*^e w_* v_*^d$ and, therefore, 
$w_*$ are elementary by definition. Thus $w'$ is the
product of $\ell$~elementary words.
\end{proof}

We can now conclude as for the length of the words
obtained using right reversing.

\begin{proof} [Proof of
Proposition~\ref{P:Bound}$(i)$]
Let $w$ be an $n$~strand braid word of
length~$\ell$. We can write it $w_1^{\epsilon_1} \pp
w_\ell^{\epsilon_\ell}$ where each~$w_k$ is a single
letter~$\ss i$ and $\epsilon_k$ is~$\pm 1$. By
Lemma~\ref{L:Decomposition}, every word~$w'$ obtained
from~$w$ by right reversing is a product of at
most~$\ell$ elementary $n$~strand words. By
Lemma~\ref{LengthElem}, each of these words has
length~$\frac 12 3^n$ at most.
\end{proof}


\section{Including monotone equivalence}

When monotone equivalence enters the picture, the
previous argument fails: simple factors may
be changed, and Lemma~\ref{L:Decomposition} does not
extend. For instance, the word~$\ss1 \ss3 \ss1
\ss3$ is a product of two simple words,
namely~$(\ss1 \ss3) (\ss1 \ss3)$, but a monotone 
equivalence transforms it to $\ss1 \ss1 \ss3 \ss3$, 
which cannot be decomposed better than 
$(\ss1) (\ss1 \ss3) (\ss3)$, a product of three simple
words. So new arguments are needed.

By Proposition~\ref{P:Basic}, every braid word~$w$ is
right reversible to a unique positive--negative
word~$v u\inv$. We shall denote~$v$ by~$\NR(w)$ (the
right numerator) and~$u$ by~$\DR(w)$ (the right
denominator). So $w \rvr \NR(w) \DR(w)\inv$ always
holds.

We first consider the case of positive--negative
words.

\begin{proof} [Proof of
Proposition~\ref{P:Bound}$(ii)$] Let $w$ be an
$n$~strand braid word of length~$\ell$. By
Proposition~\ref{P:Basic}, $w$ is right reversible
to~$\NR(w) \DR(w)\inv$, and the latter word has
length at most $(\frac12 n(n-1)-1)\ell$. In order to
prove the expected result, it is enough to prove the
following: if~$w'$ is any positive--negative word
(\ie, a word satisfying  $w'=N_R(w') D_R(w')^{-1}$)
which can be obtained from~$w$ by right reversing and
monotone equivalence, then $\NR(w')$ is equivalent
to~$\NR(w)$, and $\DR(w')$ is equivalent to~$\DR(w)$;
so, in particular, they have the same length. For an
induction, it is enough to assume that only one
monotone equivalence is used in the transformation
of~$w$ into~$w'$.

Let us display this monotone equivalence. The
hypothesis is that there exist words~$w_1, w_2$
and equivalent positive words~$v_0, v'_0$
satisfying $w \rvr w_1 v_0 w_2$, and $w_1 v'_0 w_2
\rvr w'$---the case when $v_0$ and $v'_0$ are
equivalent negative words would be treated
similarly. Proposition~\ref{P:Basic} states in
particular that the order of reversing steps does not
matter for the positive--negative word finally
obtained, so $w_1 v_0 w_2 \rvr \NR(w) \DR(w)\inv$
holds. Let us compare the reversing processes
from~$w_1 v_0 w_2$ to $\NR(w) \DR(w)\inv$ and that
from~$w_1 v'_0 w_2$ to~$w'$.

\begin{figure}[htb]
\begin{picture}(84,23)(0, 0)
\put(0,2.5){\includegraphics{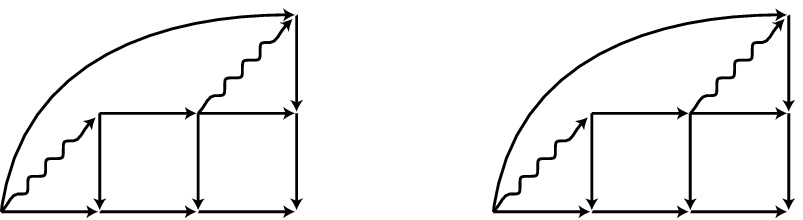}}
\put(4,0){$v_1$}
\put(14,0){$v_3$}
\put(24,0){$v_4$}
\put(5.5,6){$w_1$}
\put(11,7){$u_1$}
\put(21,7){$u_3$}
\put(31,17){$u_2$}
\put(31,7){$u_4$}
\put(25.5,16){$w_2$}
\put(14,15){$v_0$}
\put(24,11){$v_2$}
\put(5,19){$w$}
\put(54,0){$v_1$}
\put(64,0){$v'_3$}
\put(74,0){$v'_4$}
\put(55.5,6){$w_1$}
\put(61,7){$u_1$}
\put(71,7){$u'_3$}
\put(81,17){$u_2$}
\put(81,7){$u'_4$}
\put(75.5,16){$w_2$}
\put(64,15){$v'_0$}
\put(74,11){$v_2$}
\put(55,19){$w$}
\end{picture}
\caption{\smaller Comparing two reversing
processes: when $v_0$ is replaced with the
equivalent word~$v'_0$, the sequel of the reversing
changes, but the words are replaced with equivalent
words}
\label{F:Comparison}
\end{figure}

Let us introduce positive words~$v_1$, $u_1$, $v_2$,
$u_2$, $v_3$, $u_3$, $v_4$, $u_4$ satisfying $w_1 \rvr
v_1 u_1\inv$, $w_2 \rvr v_2 u_2\inv$, $u_1\inv v_0
\rvr v_3 u_3\inv$, and, finally, $u_3\inv v_2 \rvr
v_4 u_4\inv$. Then, by construction, we have $\NR(w)
= v_1 v_3 v_4$ and
$\DR(w) = u_2 u_4$ (Figure~\ref{F:Comparison}). When
we replace $v_0$ with~$v'_0$, we obtain new positive
words~$u'_3, v'_3$ satisfying $u_1\inv v'_0 \rvr
v'_3 {u'_3}\inv$, and Lemma~\ref{L:Lcm} guarantees
that $u'_3$ is equivalent to~$u_3$ and $v'_3$ is
equivalent to~$v_3$. Then, we have ${u'_3}\inv v_2
\rvr v'_4 {u'_4}\inv$ for some~$u'_4, v'_4$, and
Lemma~\ref{L:Lcm} guarantees that $u'_4$ is
equivalent to~$u_4$ and $v'_4$ is equivalent
to~$v_4$. We conclude that $\NR(w')$, which is $v_1
v'_3 v'_4$, is equivalent to~$v_1 v_3 v_4$, \ie,
to~$\NR(w)$, and that $\DR(w')$, which is~$u_2 u'_4$,
is equivalent to~$u_2 u_4$, \ie, to~$\DR(w)$.
\end{proof}

Once again, the previous argument only deals
with the final, positive--negative words
obtained using right reversing and monotone
equivalence, and it says nothing about the length
of the intermediate words. In order to prove 
Proposition~\ref{P:Bound}~$(iii)$, we need
a new argument. In the sequel, we use~$\NL(w)$
and~$\DL(w)$ for the unique positive words satisfying
$w \rvl \DL(w)\inv \NL(w)$ (the left numerator and
denominator). We recall that a word is called simple
if it is positive and represents a divisor of~$\D_n$.

\begin{lemm} \label{L:Decomposition2}
Let $w$ be a word containing $p$~negative letters and
$q$~positive letters. Then $\NR(w)$ and $\NL(w)$ are
the products of at most~$q$ simple words, and $\DR(w)$
and $\DL(w)$ are the products of at most~$p$ simple
words.
\end{lemm}

\begin{proof}
As in the proof of Lemma~\ref{L:Decomposition},
consider the right reversing grid of~$w$. It
has height~$p$ and width~$q$ and all arrows
wear simple labels. So do in particular the bottom and
right sides. This means that $\NR(w)$ is the product
of at most $q$~simple words, and, similarly,
$\DR(w)$ is the product of at most $p$~simple words.
The argument is symmetric for left reversing.
\end{proof}

\begin{proof} [Proof of
Proposition~\ref{P:Bound}$(iii)$] Let~$w$ be an
arbitrary $n$~strand braid word of length~$\ell$.
Then $w$ is right reversible to~$\NR(w) \DR(w)\inv$,
and, symmetrically, it is left reversible
to~$\DL(w)\inv \NL(w)$. For~$u$ a braid word, we
shall denote by~$\cl{u}$ the braid represented by~$u$.
Now we define $\Gamma(w)$ to be the restriction
of the Cayley graph of~$B_n$ to the divisors
of~$\cl{\DL(w) \NR(w)}$ in~$B_n^+$, \ie, 
$\Gamma(w)$ is a finite graph, containing precisely
those vertices that lie on some geodesic path from~$1$
to~$\cl{\DL(w) \NR(w)}$---these paths all have the
same length, since they correspond to positive words
equivalent  to~$\DL(w) \NR(w)$. By
Lemma~\ref{L:Decomposition2}, the word~$\DL(w)
\NR(w)$ is the product of at most~$\ell$
simple words, hence its length is at most $\frac12
n(n-1) \ell$.

Let~$\beta$ be a vertex of the graph~$\Gamma(w)$. We
say that a braid word~$u$ is {\it traced}
in~$\Gamma(w)$ from~$\beta$ if there exists a path
labelled~$u$ starting at~$\beta$ in~$\Gamma(w)$,
\ie, we can read all letters of~$u$ successively
without leaving~$\Gamma(w)$. Then it is proved
in~\cite{Dfo} that the word~$w$ itself is traced
from~$\cl{\DL(w)}$ in~$\Gamma(w)$, and that the
family of all words traced from a fixed vertex
in~$\Gamma(w)$ is closed under right and left
reversing, and under monotone equivalence.
Therefore, every word~$w'$ that can be derived
from~$w$ using reversing and monotone equivalence is
traced from~$\cl{\DL(w)}$ in~$\Gamma(w)$.

Now, we attribute a weight to every edge~$e$
in~$\Gamma(w)$, namely the integer~$F_d$, where
$d$ is the distance from the source vertex of~$e$ to
the final vertex of~$\Gamma(w)$, and $F_d$ is the
$d$th Fibonacci number: $F_1 = F_2 = 1$, and $F_k =
F_{k-1} + F_{k-2}$ for $k \ge 2$
(Figure~\ref{F:Weight}). Finally we define the weight
of a path in~$\Gamma(w)$ to be the sum of the weight
of its edges.

\begin{figure}[htb]
\begin{picture}(90,28)(0, 0)
\put(0,0){\includegraphics{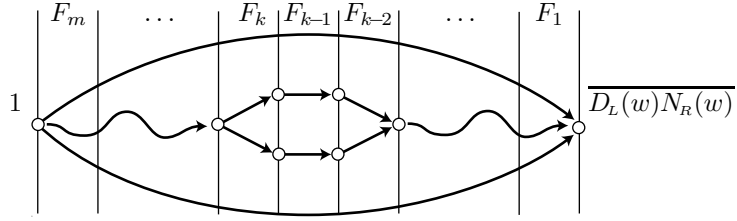}}
\put(2.5,26){$F_m$}
\put(15,26){$\pp$}
\put(55,26){$\pp$}
\put(27.5,26){$F_k$}
\put(33.5,26){$F_{k\!-\!1}$}
\put(41.5,26){$F_{k\!-\!2}$}
\put(67,26){$F_1$}
\put(-3,14){$1$}
\put(74,14){$\cl{\DL(w) \NR(w)}$}
\end{picture}
\caption{\smaller Weights for the edges in the
Cayley graph~$\Gamma(w)$}
\label{F:Weight}
\end{figure}

Then we observe that the weight cannot increase when
right reversing or monotone equivalence is
performed. Indeed, reversing $\sss i \ss j$ to $\ss j 
\ss i \sss j \sss i$ (with $\vert i - j\vert = 1$)
replaces two edges contributing say $2F_k$ to the
total weight with four edges contributing $2F_{k-1}
+ 2F_{k-2}$, \ie, $2F_k$ again. Similarly, reversing
$\sss i \ss j$ to $\ss j \sss i$ (with $\vert i -
j\vert \ge 2$) diminishes the contribution to the
weight from~$2F_k$ to~$2F_{k-1}$, and deleting 
$\sss i \ss i$ diminishes it by~$2F_k$. Finally,
replacing some positive (\resp negative) subword with
an equivalent positive (\resp negative) word
preserves the weight. As each letter in a braid word
contributes at least~$1$ in the weight, we deduce
that the length of any word obtained from~$w$ using
right reversing and monotone equivalence is bounded
above by the weight of~$w$.

The latter is the sum of~$\ell$ Fibonacci numbers
between~$F_1$ and~$F_m$, where $m$ is at most
$\frac12 n(n-1) \ell$. One easily checks that the
worst case is when $w$ consists of $\ell/2$ negative
letters and $\ell/2$ positive letters with weights
$F_m$, $F_{m-1}$, \pp, $F_{m - \ell/2 + 1}$.
Using the very rough estimate $F_k \le 2^{k-2}$ for
$k\ge 2$ we obtain an upper bound of
$$
2\sum_{k=1}^{m}F_k \le 2(1+\sum_{k=0}^{m-2} 2^k) = 2^m
\le 2^{\frac12 n(n-1) \ell}
$$
on the weight of $w$.
\end{proof}


\section{A new conjecture}

The previous results leave open all questions 
about simultaneous left and right reversing. We 
shall conclude with a conjecture which does appear 
to be true---it is confirmed by extensive computer
experiments---and which would still imply a linear
upper bound on the space complexity of the handle
reduction algorithm.

First we recall that reduced positive equivalences
can be  decomposed into commutation relations---like
$\sigma_1 \sigma_3 \to \sigma_3 \sigma_1$---and
Reidemeister~III relations---for instance
$\sigma_2\sigma_3\sigma_2 \to \sigma_3 \sigma_2
\sigma_3$---each followed by a free reduction.

\begin{conj} \label{C:NewConj}
Let $w$ be an $n$~strand braid word of length~$\ell$.
Let $w'$ be another such word which is obtained from
$w$ by a sequence of reduced  commutation relations
and reduced word reversings.  Then $w'$ is of length
at most $C_n \ell$, where 
$C_n$ is a constant depending only on $n$.
\end{conj}

We have no good guess, however, what the constant
$C_n$ should be. We do know that it is not
$\frac12 n(n-1)$, as the $4$~strand braid word
$\s_3^2 \sss2 \s_1^2 \sss2 \ss3 \ss1$ of length 8
can be transformed into a word of length 52 (which is
larger than $6\cdot 8 = 48$).


We also remark that there is a slightly weaker
version of Conjecture~\ref{C:NewConj} where
``reduced'' is replaced with ``strongly reduced''
everywhere.

In order to prove Conjecture~\ref{C:NewConj}, it might
be useful to consider Bestvina's product structure on
the flag complex $\widehat{\mathcal{X}}_D$ which is
closely related to the Cayley graph of $B_n$
\cite{Bes,CMW}. Bestvina showed that there is a
natural homeomorphism $\widehat{\mathcal{X}}_D \cong 
\mathcal{X}_D\times\mathbb{R}$,
where $\mathcal{X}_D$ is another complex which
satisfies a certain weak non-positive curvature
condition, and is conjectured to be $C\!AT(0)$. We
think now of a braid word as a path in
$\widehat{\mathcal{X}}_D$, and of our transformations
of braid words as deformations of the path that
preserve its endpoints. Then banning positive
equivalences of Reidemeister~III type amounts to
forbidding the most obvious way of deforming a path
in the
$\mathcal{X}_D$-direction. In other words, applying
only commutation relations and word reversings
means deforming the path mainly in the $\mathbb
R$-direction.


\end{document}